\documentclass[reqno, 11pt]{smfart}
\usepackage{srcltx}
\usepackage[utf8]{inputenc}

\usepackage{a4wide}
\usepackage{amsfonts, amsthm, amsmath, amssymb}
\usepackage[francais,english]{babel}
\usepackage{amsmath, amsfonts, amssymb, amsthm}
\usepackage{graphicx, subfigure}
\usepackage{fancyhdr}
\usepackage{stmaryrd}
\usepackage{enumerate}
\usepackage{shorttoc}
\usepackage{xy}
\usepackage{variations}
\def\BbN{\mathbb N}
\newcommand{\lb}{\label}
  \newcommand{\ds}{\displaystyle}

\theoremstyle{plain}

\renewcommand{\theequation}{\thesection.\arabic{equation}}

\author{Abdelmajid Ben Hadj Salem}
\address{6, rue du Nil, Cit\'e Soliman Er-Riadh\\ 8020 Soliman\\ Tunisie } 
\email{abenhadjsalem@gmail.com}

\def\si {\sigma}
\def\al {\alpha}
\def\BbR{{\mathbb R}}

\def\lb {\label}
\def\BbN{{\mathbb N}}
 
\def\BbC{{\mathbb C}}
\def\be{\begin{equation}}
\def\ee{\end{equation}}
\def\ba{\begin{eqnarray}}
\def\ea{\end{eqnarray}}
\def\ds{\displaystyle}
\def\bq{\begin{quote}}
\def\eq{\end{quote}}
\usepackage{graphicx}
\newtheorem{theorem}{Th\'eor\`eme}[section]
\newtheorem{conjecture}[theorem]{\textbf{Conjecture}}
\newtheorem{proposition}[theorem]{\textbf{Proposition}}
\newtheorem{equivalence}[theorem]{\textbf{Equivalence}}

\begin{document}
\renewcommand{\contentsname}{Sommaire}
\renewcommand{\abstractname}{Abstract}
\setcounter{section}{0}
\rhead{}
\chead{}
\lhead{}
\fancyhead[CO]{\footnotesize{ \textsc{Une Solution de l'Hypoth\`ese de Riemann }}}
\fancyhead[CE]{\footnotesize{\textsc{A. Ben Hadj Salem}}}
\altkeywords{Zeta function; nontrivial zeros of zeta function; equivalent statements; definition of limit of real sequences.}
\title{Une Solution de l'Hypoth\`ese de Riemann}
\alttitle{A Solution of Riemann Hypothesis}
	\dedicatory{A mon \'epouse Wahida}
\maketitle


%
%
%

\selectlanguage{french}
\begin{abstract} 
En 1859, Georg Friedrich Bernhard Riemann avait annonc\'e la conjecture suivante \cite{bib1}, dite Hypoth\`ese de Riemann: \textit{Les z\'eros non triviaux $s=\sigma+it$ de la fonction zeta d\'efinie par:  
$$\zeta(s) = \sum_{n=1}^{+\infty}\frac{1}{n^s},\,\mbox{pour }\quad \Re(s)>1$$
ont comme parties r\'eelles $\sigma= \ds \frac{1}{2}$}.

On donne une d\'emonstration que $\sigma= \ds \frac{1}{2}$ en utilisant une proposition \'equivalente de l'Hypoth\`ese de Riemann.
     \end{abstract}
\begin{altabstract}
In 1859, Georg Friedrich Bernhard Riemann had announced the following conjecture \cite{bib1}, called Riemann Hypothesis : \textit{The nontrivial roots (zeros) $s=\sigma+it$ of the zeta function, defined by:  
$$\zeta(s) = \sum_{n=1}^{+\infty}\frac{1}{n^s},\,\mbox{for}\quad \Re(s)>1$$
have real part $\sigma= \ds \frac{1}{2}$}.

We give a proof that $\sigma= \ds \frac{1}{2}$ using an equivalent statement of Riemann Hypothesis.
\end{altabstract}


\vspace*{6pt}\tableofcontents  
\mainmatter
\section{Introduction}
En 1859, G.F.B. Riemann avait annonc\'e la conjecture suivante \cite{bib1}: 
\begin{conjecture}
\textit{Soit $\zeta(s)$ la fonction complexe de la variable complexe $s=\sigma+it$ d\'efinie par le prolongement analytique de la fonction: 
$$ \zeta_1(s)=\sum_{n=1}^{+\infty}\frac{1}{n^s},\,\mbox{pour}\,\,\Re(s)=\sigma>1 $$ 
sur tout le plan complexe sauf au point $s=1$. Alors les z\'eros non triviaux de $\zeta(s)=0$ sont de la forme:}
$$ s=\frac{1}{2}+it $$ 
\end{conjecture}
Dans cette communication, nous donnons une d\'emonstration que $ \sigma= \ds \frac{1}{2}$. Notre id\'ee est de partir  d'une proposition \'equivalente de l'Hypoth\`ese de Riemann et en utilisant la d\'efinition de la limite des suites r\'eelles.
\
\subsection{La fonction $\zeta$}
Notons par $s=\sigma+it$ la variable complexe de $\BbC$. Pour $\Re(s)=\sigma>1$, appelons $\zeta_1$ la fonction d\'efinie par :
$$ \zeta_1(s)=\sum_{n=1}^{+\infty}\frac{1}{n^s},\,\mbox{avec}\,\,\Re(s)=\sigma>1 $$ 
Nous savons qu'avec la d\'efinition pr\'ec\'edente, la fonction $\zeta_1$ est une fonction analytique de $s$. Notons par $\zeta(s)$ la fonction obtenue par prolongement analytique de $\zeta_1(s)$, alors nous rappelons le th\'eor\`eme suivant \cite{bib2}:
\begin{theorem} \textit{Les z\'eros de $\zeta(s)$ satisfont :}

1. \textit{$\zeta(s)$ n'a pas de z\'eros pour $\Re(s)>1$;}

2. \textit{le seul p\^ole de $\zeta(s)$ est au point $s=1$; son r\'esidu vaut 1 et il est simple;}

3. \textit{les z\'eros triviaux de $\zeta(s)$ sont d\'etermin\'es pour les valeurs $s=-2,-4,\ldots$;}

4. \textit{les z\'eros non triviaux se situent dans la r\'egion $0\leq \Re(s)\leq1$ dite bande critique et ils sont sym\'etriques respectivement par rapport \`a l'axe vertical $\Re(s)=\ds \frac{1}{2}$ et l'axe des r\'eels $\Im(s)=0$.}\lb{th2} 
\end{theorem}
Par suite, la conjecture relative \`a l'Hypoth\`ese de Riemann est exprim\'ee comme suit:
\begin{conjecture}\textit{(Hypoth\`ese de Riemann,\cite{bib2}) Tous les z\'eros non triviaux de  $\zeta(s)$ sont sur la droite critique $\Re(s)=\ds \frac{1}{2}$. }
\end{conjecture}
En plus des propri\'etes cit\'ees par le th\'eor\`eme cit\'e  ci-dessus, la fonction $\zeta(s)$ v\'erifie la relation fonctionnelle \cite{bib2} pour $s\neq1$:
\be
\zeta(1-s)=2^{1-s}\pi^{-s}cos\ds \frac{s\pi}{2}\Gamma(s)\zeta(s) \lb{rim5}
\ee
o\`u $\Gamma(s)$ est la fonction d\'efinie sur le demi-plan $\Re(s)>0$ par:
$$ \Gamma(s)=\int_0^{\infty} e^{-t}t^{s-1}dt,\quad  $$ 
Alors, au lieu d'utiliser la fonctionnelle donn\'ee par (\ref{rim5}), nous allons utiliser celle pr\'esent\'ee par G.H. Hardy \cite{bib3} \`a  savoir  la fonction eta de Dirichlet \cite{bib2}:
$$ \eta(s)=\ds \sum_{n=1}^{+\infty}\frac{(-1)^{n-1}}{n^s}=(1-2^{1-s})\zeta(s) $$ 
Elle est convergente pour tout $s\in \BbC$ avec $\Re(s)>0$.
\subsection{Une Proposition \'equivalente \`a l'Hypoth\`ese de Riemann}
Parmi les propositions \'equivalentes \`a l'Hypoth\`ese de Riemann celle de la fonction eta de Dirichlet qui s'\'ennonce comme suit \cite{bib2}:
\begin{equivalence}
\textit{L'Hypoth\`ese de Riemann est \'equivalente \`a l'\'ennonc\'e que tous les z\'eros de la fonction eta de Dirichlet:
$$ \eta(s)=\ds \sum_{n=1}^{+\infty}\frac{(-1)^{n-1}}{n^s}=(1-2^{1-s})\zeta(s) $$ 
qui se situent dans la bande critique $0 <\Re(s) < 1$, sont sur la droite critique $\Re(s)=\ds \frac{1}{2}$.}
\end{equivalence}\lb{eq1}

\section{D\'emonstration  que les z\'eros de la $\eta(s)$ sont sur la droite critique $ \Re(s)=\ds \frac{1}{2}$ }
Notons par $s=\sigma+it$ avec $0<\sigma<1$. Consid\'erons maintenant un z\'ero de $\eta(s)$ qui se trouve dans la bande critique et appelons $s=\sigma+it$ ce z\'ero, nous avons donc $0<\sigma<1$ et $\eta(s)=0\Longrightarrow (1-2^{1-s})\zeta(s)=0$. Notons $\zeta(s)=A+iB$, et $\theta=t Log2$,  alors:
$$(1-2^{1-s})\zeta(s)=\left[A(1-2^{1-\sigma}cos\theta)-2^{1-\sigma}Bsin\theta\right]+i\left[B(1-2^{1-\sigma}cos\theta)+2^{1-\sigma}Asin\theta\right]$$
$(1-2^{1-s})\zeta(s)=0$ donne le syst\`eme:
\ba
A(1-2^{1-\sigma}cos\theta)-2^{1-\sigma}Bsin\theta=0 \nonumber \\
B(1-2^{1-\sigma}cos\theta)+2^{1-\sigma}Asin\theta=0 \nonumber
\ea
Comme les fonctions \textit{sin} et \textit{cos} ne s'annulent pas simultan\'ement, supposons par exemple que $sin\theta \neq0$, la premi\`ere \'equation du syst\`eme donne $B=\ds \frac{A(1-2^{1-\sigma}cos\theta)}{2^{1-\sigma}sin\theta}$, la deuxi\`eme \'equation s'\'ecrit:
$$ \frac{A(1-2^{1-\sigma}cos\theta)}{2^{1-\sigma}sin\theta}(1-2^{1-\sigma}cos\theta)+2^{1-\sigma}Asin\theta=0\Longrightarrow A=0$$
Par suite, $B=0\Longrightarrow \zeta(s)=0$, il s'ensuit que:
\be
  \fbox{$s\,\mbox{\textit{est un z\'ero de}}\, \eta(s)\,\,\mbox{\textit{dans la bande critique est aussi un z\'ero de}}\, \zeta(s)$} \lb{rim100}
	\ee
Reciproquement, si $s$ est un z\'ero de $\zeta(s)$ dans la bande critique, soit $\zeta(s)=A+iB=0\Longrightarrow \eta(s)=(1-2^{1-s})\zeta(s)=0$, donc $s$ est aussi un z\'ero de $\eta(s)$ dans la bande critique. Nous pouvons \'ecrire:
\be
  \fbox{$s\,\mbox{\textit{est un z\'ero de}}\, \zeta(s)\,\,\mbox{\textit{dans la bande critique est aussi un z\'ero de}}\, \eta(s)$} \lb{rim100a}
	\ee

Ecrivons la fonction $\eta$:
\ba
&\eta(s)=\ds \sum_{n=1}^{+\infty} \frac{(-1)^{n-1}}{n^s}=\sum_{n=1}^{+\infty}(-1)^{n-1} e^{-sLogn}=\sum_{n=1}^{+\infty} (-1)^{n-1}e^{-(\sigma+it)Logn}=\nonumber &\\ & \ds = \sum_{n=1}^{+\infty} (-1)^{n-1}e^{-\sigma Logn}. e^{-it Logn} & \nonumber \\ & \ds  =\sum_{n=1}^{+\infty} (-1)^{n-1}e^{-\sigma Logn} (cos(t Logn) - isin( t Logn)) &\nonumber 
\ea 
D\'efinissons la suite de fonctions $((\eta_n)_{n\in \BbN^*}(s))$, par:
$$\eta_n(s)=\ds \sum_{k=1}^{n} \frac{(-1)^{k-1}}{k^s}=\sum_{k=1}^{n}(-1)^{k-1} \frac{cos(t Logk)}{k^{\sigma}} - i\sum_{k=1}^n(-1)^{k-1}\frac{sin( t Logk)}{k^{\sigma}} $$
avec $s=\sigma+it$ et $t\neq0$. 
\\

Soit $s$ un z\'ero de $\eta$ dans la bande critique, soit $\eta(s)=0$, avec $0<\sigma<1$. Par suite, on peut \'ecrire $lim_{n\longrightarrow +\infty}\eta_n(s)=0=\eta(s)$. Ce qui donne:
\ba
lim_{n\longrightarrow +\infty}\ds \sum_{k=1}^{n} (-1)^{k-1}\frac{cos(t Logk)}{k^{\sigma}}=0 \nonumber \\
lim_{n\longrightarrow +\infty}\ds \sum_{k=1}^{n} (-1)^{k-1}\frac{sin(t Logk)}{k^{\sigma}}=0 \nonumber 
\ea
Utilisons la d\'efinition de la limite d'une suite, on peut \'ecrire:
\ba
\forall \epsilon_1>0 \quad \exists n_r,\forall N>n_r \quad |\Re(\eta(s)_N)|<\epsilon_1 \nonumber \\
\forall \epsilon_2>0 \quad \exists n_i,\forall N>n_i \quad |\Im(\eta(s)_N)|<\epsilon_2 \nonumber 
\ea
En prenant $\epsilon=\epsilon_1=\epsilon_2$ et $N>max(n_r,n_i)$, on obtient:
\ba
0<\ds \sum_{k=1}^{N} \frac{cos^2(t Logk)}{k^{2\sigma}}+2\sum_{k,k'=1;k<k'}^{N} \frac{(-1)^{k+k'}cos(t Logk).cos(t Logk')}{k^{\sigma}k'^{\sigma}}<\epsilon^2 \nonumber  \\
0<\ds \sum_{k=1}^{N} \frac{sin^2(t Logk)}{k^{2\sigma}}+2\sum_{k,k'=1;k<k'}^{N} \frac{(-1)^{k+k'}sin(t Logk).sin(t Logk')}{k^{\sigma}k'^{\sigma}}<\epsilon^2 \nonumber 
\ea
En faisant la somme des deux derni\`eres in\'egalit\'es, on obtient:
\be
0<\sum_{k=1}^{N}\ds \frac{1}{k^{2\sigma}}+2\sum_{k,k'=1;k< k'}^{N} (-1)^{k+k'}\frac{cos(t Log(k/k'))}{k^{\sigma}k'^{\sigma}}<2\epsilon^2 \lb{rim16} 
\ee
\subsection{Cas $\sigma=\ds \frac{1}{2}\Longrightarrow 2\sigma=1$}
On suppose que $\si=\ds \frac{1}{2} \Longrightarrow 2\si=1$. Commen\c{c}ons par rappeler le th\'eor\`eme de Hardy (1914) \cite{bib2},\cite{bib3}:
\begin{theorem}
\textit{Il y'a une infinit\'e de z\'eros de $\zeta(s)$ sur la droite critique.}
\end{theorem}
Des propositions (\ref{rim100}-\ref{rim100a}), nous d\'eduisons la proposition suivante:
\begin{proposition}
\mbox{Il y'a une infinit\'e de z\'eros de} $\eta(s)$ \mbox{sur la droite critique.}
\end{proposition}
Soit $s_j=\frac{1}{2}+it_j$ un des z\'eros de la fonction $\eta(s)$ sur la droite critique, soit $\eta(s_j)=0$. L'\'equation (\ref{rim16}) s'\'ecrit pour $s_j$:
$$0<\sum_{k=1}^{N}\ds \frac{1}{k}+2\sum_{k,k'=1;k<k'}^{N} (-1)^{k+k'}\frac{cos(t_j Log(k/k'))}{\sqrt{k}\sqrt{k'}}<2\epsilon^2 $$ 
ou encore:
$$\sum_{k=1}^{N}\ds \frac{1}{k}<2\epsilon^2-2\sum_{k,k'=1;k<k'}^{N} (-1)^{k+k'}\frac{cos(t_j Log(k/k'))}{\sqrt{k}\sqrt{k'}} $$ 
Si on fait tendre $N$ vers $+\infty$, la s\'erie $\ds \sum_{k=1}^{N} \frac{1}{k}$ est divergente et devient infinie. Soit:
$$ \sum_{k=1}^{+\infty}\ds \frac{1}{k}\leq 2\epsilon^2-2\sum_{k,k'=1;k<k'}^{+\infty} (-1)^{k+k'}\frac{cos(t_j Log(k/k'))}{\sqrt{k}\sqrt{k'}} $$
 Par suite, nous obtenons le r\'esultat important suivant:
\be
 \fbox{ $lim_{N\longrightarrow +\infty} \ds \sum_{k,k'=1;k<k'}^{N}\ds(-1)^{k+k'}  \frac{cos(t_j Log(k/k'))}{\sqrt{k}\sqrt{k'}}=-\infty $} \lb{rim19a}
\ee 
sinon, nous aurons une contradiction avec le fait que : 
$$ lim_{N\longrightarrow +\infty} \ds \sum_{k=1}^{N}\ds (-1)^{k-1}\frac{1}{k^{s_j}}=0$$
Comme $t_j\neq0$, et qu'il y'a une infinit\'e de z\'eros sur la droite critique, alors le r\'esultat de la formule donn\'ee par (\ref{rim19a}) est ind\'ependant de $t_j$. Revenons maintenant \`a $s=\sigma+it$ un z\'ero de $\eta(s)$ dans la bande critique, soit $\eta(s)=0$. Prenons $\sigma=\ds \frac{1}{2}$. En partant de la d\'efinition de la limite des suites appliqu\'ee ci-dessus, nous obtenons:
$$\sum_{k=1}^{+\infty}\ds \frac{1}{k}\leq 2\epsilon^2-2\sum_{k,k'=1;k<k'}^{+\infty} (-1)^{k+k'}\frac{cos(t Log(k/k'))}{\sqrt{k}\sqrt{k'}}$$ 
avec sans aucune contradiction. De la proposition (\ref{rim100}) il s'ensuit que $\zeta(s)=\zeta(\frac{1}{2}+it)=0$. Il existe donc des z\'eros de $\zeta(s)$ sur la droite critique $\Re(s)=\ds \frac{1}{2}$.
 
\subsection{Cas $\ds 0<\si<\frac{1}{2}$}
\subsubsection{Cas o\`u il n'existe pas de z\'eros de $\eta(s)$ avec $s=\sigma+it$ et $\ds 0<\si<\frac{1}{2}$}
En utilisant, pour ce cas, le point 4 du th\'eor\`eme (\ref{th2}), nous d\'eduisons que la fonction $\eta(s)$ n'a pas de z\'eros avec $s=\sigma+it$ et $\ds \frac{1}{2}<\sigma<1$. Par suite, d'apr\`es la proposition (\ref{rim100}), il s'ensuit que la fonction $\zeta(s)$ a  ses z\'eros seulement sur la droite critique $\Re(s)=\si=\ds \frac{1}{2}$ et \textbf{l'Hypoth\`ese de Riemann est vraie}. 
\subsubsection{Cas o\`u il existe des z\'eros de $\eta(s)$ avec $s=\sigma+it$ et $ 0<\si< \ds \frac{1}{2}$}
Supposons qu'il existe $s=\sigma+it$ un z\'ero de $\eta(s)$ soit $\eta(s)=0$ avec $0<\si<\frac{1}{2}\Longrightarrow s\in$ \`a la bande critique. Nous \'ecrivons l'\'equation (\ref{rim16}),:
$$0<\sum_{k=1}^{N}\ds \frac{1}{k^{2\sigma}}+2\sum_{k,k'=1;k< k'}^{N} (-1)^{k+k'}\frac{cos(t Log(k/k'))}{k^{\sigma}k'^{\sigma}}<2\epsilon^2 $$
ou:
$$ \sum_{k=1}^{N}\ds \frac{1}{k^{2\sigma}}<2\epsilon^2-2\sum_{k,k'=1;k< k'}^{N} (-1)^{k+k'}\frac{cos(t Log(k/k'))}{k^{\sigma}k'^{\sigma}} $$
Or $2\sigma<1$, il s'ensuit que $lim_{N\longrightarrow+\infty} \ds \sum_{k=1}^{N} \frac{1}{k^{2\sigma}}$ tende vers $+\infty$ et nous obtenons:
$$ \sum_{k,k'=1;k<k'}^{+\infty}\ds(-1)^{k+k'}  \frac{cos(t Log(k/k'))}{k^{\sigma}k'^{\sigma}}=-\infty $$
\subsection{Cas $\ds \frac{1}{2}<\Re(s)<1$}\lb{23}
Soit $s=\si+it$ le z\'ero de $\eta(s)$ dans $0<\Re(s)<\ds \frac{1}{2}$, objet du paragraphe pr\'ec\'edent. Suivant le point 4 du th\'eor\`eme \ref{th2}, le nombre complexe $s'=1-\si+it$ est aussi un z\'ero de la fonction $\eta(s)$ dans la bande $\ds \frac{1}{2}<\Re(s)<1$. En appliquant (\ref{rim16}), nous avons:  
\be
0<\sum_{k=1}^{N}\ds \frac{1}{k^{2\si'}}+2\sum_{k,k'=1;k< k'}^{N} (-1)^{k+k'}\frac{cos(t Log(k/k'))}{k^{\si'}k'^{\si'}}<2\epsilon^2 \lb{rim16b} 
\ee
Or $2\si'=2(1-\si)>1\Longrightarrow \si<\ds \frac{1}{2}$, la s\'erie $\sum_{k=1}^{N}\ds \frac{1}{k^{2\si'}}$ est convergente vers une constante $C(\si')$. De l'\'equation (\ref{rim16b}), nous d\'eduisons que :
$$\sum_{k,k'=1;k< k'}^{+\infty} (-1)^{k+k'}\frac{cos(t Log(k/k'))}{k^{\si'}k'^{\si'}}=-\frac{C(\si')}{2}>-\infty $$ 
Maintenant fixons $t=\Im(s')$ et consid\'erons la fonction $F_N(u)$ d\'efinie par:
\ba
F_N(u)=\sum_{k,k'=1;k< k'}^{N} (-1)^{k+k'}\frac{cos(t Log(k/k'))}{k^{u}k'^{u}}=\nonumber \\
= \sum_{k,k'=1;k< k'}^{N} (-1)^{k+k'}cos(t Log(k/k'))e^{-uLog(kk')},\quad u\in]0,1[ \nonumber 
\ea
La fonction $F_N(u)$ est continue pour $\forall N\in \BbN^*$ et $u\in]0,1[$, et nous avons obtenu pr\'ec\'edement que pour $N\longrightarrow +\infty$  :
$$ \left\{\begin{array}{l}
\ds \sum_{k,k'=1;k< k'}^{+\infty} (-1)^{k+k'}\frac{cos(t Log(k/k'))}{k^{\si'}k'^{\si'}}=-\frac{C(\si')}{2}\quad \mbox{pour}\,\,u=\si'=1-\si >\ds \frac{1}{2}\\ \\
\ds \sum_{k,k'=1;k< k'}^{+\infty} (-1)^{k+k'}\frac{cos(t Log(k/k'))}{\sqrt{k}\sqrt{k'}}=-\infty\quad \mbox{pour }\,\,u=\frac{1}{2} \\
\\
\ds \sum_{k,k'=1;k< k'}^{+\infty} (-1)^{k+k'}\frac{cos(t Log(k/k'))}{k^{\si}k'^{\si}}=-\infty \quad \mbox{pour}\,\,u=\si < \ds \frac{1}{2}
\end{array}\right. $$
Ecrivons que $F_N(u)$ est continue au point $u=1/2$, on peut \'ecrire:
$$\forall \epsilon>0, \exists \delta\, \mbox{tel que}\,\,\forall \,u\,   /\, |u-1/2|<\delta \Longrightarrow |F_N(u)-F_N(1/2)|<\epsilon$$
Soit $u=\sigma'\in ]0,1[$ avec $\si'>\frac{1}{2}$ v\'erifiant $\sigma'-\frac{1}{2}<\delta$, on a alors l'\'equation:
\ba
&|F_N(\sigma)-F_N(1/2)|<\epsilon\Longrightarrow - \epsilon+F_N(\sigma')<\ds \sum_{k,k'=1;k< k'}^{N} (-1)^{k+k'}\frac{cos(t Log(k/k'))}{\sqrt{k}\sqrt{k'}}<\epsilon +F_N(\sigma') \nonumber &\\&
\Longrightarrow -\epsilon+\ds \sum_{k,k'=1;k< k'}^{N} (-1)^{k+k'}\frac{cos(t Log(k/k'))}{k^{\si'}k'^{\si'}}<\ds \sum_{k,k'=1;k< k'}^{N} (-1)^{k+k'}\frac{cos(t Log(k/k'))}{\sqrt{k}\sqrt{k'}} & \nonumber
\ea
Comme pour $t,u$ fix\'es, la fonction $F_N$  est d\'efinie pour tout entier $N>0$, faisons alors tendre $N$ vers $+\infty$, nous obtenons:
\ba
& -\epsilon+\ds \sum_{k,k'=1;k< k'}^{+\infty} (-1)^{k+k'}\frac{cos(t Log(k/k'))}{k^{\si'}k'^{\si'}}\leq \ds \sum_{k,k'=1;k< k'}^{+\infty} (-1)^{k+k'}\frac{cos(t Log(k/k'))}{\sqrt{k}\sqrt{k'}} & \nonumber \\
&\Longrightarrow -\epsilon \ds -\frac{C(\si')}{2}\leq \ds \sum_{k,k'=1;k< k'}^{+\infty} (-1)^{k+k'}\frac{cos(t Log(k/k'))}{\sqrt{k}\sqrt{k'}}=-\infty& \nonumber
\ea
D'o\`u la contradiction avec $C(\si')$ born\'ee. Par suite, l'hypoth\`ese qu'il existe des z\'eros de $\eta(s)$ dans l'intervalle $\ds \frac{1}{2}<\Re(s)<1$ \'etudi\'ee dans [\ref{23}] est fausse. Il s'ensuit que la fonction $\eta(s)$ ne s'annule pas dans les intervalles $\ds 0<\Re(s)<\frac{1}{2}$ et $\ds \frac{1}{2}<\Re(s)<1$ et par suite la fonction $\eta(s)$ a ses z\'eros non triviaux sur la droite critique $\Re(s)=\ds \frac{1}{2}$ de la bande critique. 

\section{Conclusion}
En r\'esum\'e: pour nos d\'emonstrations, nous avons fait usage de la fonction $\eta(s)$ de Dirichlet:
 $$\eta(s)=\ds \sum_{n=1}^{+\infty}\frac{(-1)^{n-1}}{n^s}=(1-2^{1-s})\zeta(s), \quad s=\sigma+it $$
dans la bande critique $0<\Re(s)<1$, en obtenant:

- $\eta(s)$ s'annule pour $0<\sigma=\Re(s)  = \ds \frac{1}{2}$;  

- $\eta(s)$ ne s'annule pas pour $0<\sigma=\Re(s)<\ds \frac{1}{2}$ et $\ds \frac{1}{2}<\sigma=\Re(s)<1$.
\\

Par suite, tous les z\'eros non triviaux de $\eta(s)$ dans la bande critique $0<\Re(s)<1$ s'annulent sur la droite critique $\Re(s)=\ds \frac{1}{2}$. En appliquant la proposition \'equivalente \`a l'Hypoth\`ese de Riemann \ref{eq1}, les z\'eros non triviaux de la fonction $\zeta(s)$ se trouvent sur la droite critique $\Re(s)=\ds \frac{1}{2}$. La d\'emonstration de l'Hypoth\`ese de Riemann est ainsi achev\'ee. Nous annon\c{c}ons donc le th\'eor\`eme important :
\begin{theorem}(\textit{\textbf{Abdelmajid Ben Hadj Salem, 2017):}}

\textit{Tous les z\'eros non triviaux de la fonction $\zeta(s)$ avec $s=\sigma+it$ se situent sur l'axe vertical $\Re(s)=\ds \frac{1}{2}$.}\lb{th2a}
\end{theorem} 

%





\begin{thebibliography}{99}
\bibitem{bib1} \textsc{Enrico Bombieri}. \textit{The Riemann Hypothesis}, pp 107-124. The Millennium Prize Problems. J. Carlson, A. Jaffe, and A. Wiles Editors. 160 pages. Published by The American Mathematical Society, Providence, RI, for The Clay Mathematics Institute, Cambridge, MA. 2006.

\bibitem{bib2} \textsc{Peter Borwein, Stephen Choi, Brendan Rooney and Andrea Weirathmueller}. \textit{The Riemann Hypothesis - A Resource for the Afficionado and Virtuoso Alike}. First Edition. CMS Books in Mathematics. Springer-Verlag New York. 533 pages. 2008.

\bibitem{bib3} \textsc{E.C. Titchmarsh, D.R. Heath-Brown}. \textit{The Theory of the Riemann Zeta-Function}. Second Edition revised by D.R. Heath-Brown. Oxford University Press, New York. 418 pages. 1986.
\end{thebibliography}
\end{document}